2000]57M25, 53A04, 49Q10, 52C15, 52C17

# On Thickness and Packing Density for Knots and Links

## Rob Kusner

ABSTRACT. We describe some problems, observations, and conjectures concerning thickness and packing density of knots and links in $\mathbb{S}^3$ and $\mathbb{R}^3$. We prove the thickness of a nontrivial knot or link in $\mathbb{S}^3$ is no more than $\frac{\pi}{4}$, the thickness of a Hopf link. We also give arguments and evidence supporting the conjecture that the packing density of thick links in $\mathbb{R}^3$ or $\mathbb{S}^3$ is generally less than $\frac{\pi}{\sqrt{12}}$, the density of the hexagonal packing of unit disks in $\mathbb{R}^2$.

## Introduction

This note describes some problems, observations, and conjectures concerning *thick* knots and links — that is, collections of solid tori or cylinders embedded as constant radius tubes around simple closed curves or proper arcs — in Riemannian 3-manifolds, especially in $\mathbb{R}^3$ and $\mathbb{S}^3$. The maximal such radius is called the *thickness* of the collection. In a compact 3-manifold like $\mathbb{S}^3$, it makes sense to maximize this over isotopy as well, leading to the first basic problem:

*How does the maximal thickness of a knot or link in the round three-sphere depend on its isotopy type, and what is the geometry of a thickness maximizer?*

For thick knots and links in $\mathbb{R}^3$, to accommodate the effect of rescaling, one considers instead *ropelength* — length of core curves divided by thickness — and the analogous problems about ropelength minimizers have partly been answered, thanks to the work of many people at this workshop (for instance, see [3] and the references therein, as well as [5] for the latest ropelength bounds).

The second basic problem arose from thinking about how to estimate length and thickness using the formula for the volume of a tube; the ratio of the tube volume to the ambient volume — the *packing density* — appears as a correction factor. This may be an interesting scale-invariant measure for thickness maximizing configurations, but it also would be good to have a general estimate (better than 1) for the packing density of any thick knot or link. To have the inequality go the right way, the ambient curvature should be non-negative, and one is led to ask:

*Can the packing density for a thick knot or link in Euclidean three-space, or (aside from some trivial situations) in the round three-sphere, ever exceed the maximal disk-packing density in the plane?*

Of course, this maximal density $\frac{\pi}{\sqrt{12}} = .9069\cdots$ is achieved by packing infinite solid cylinders in $\mathbb{R}^3$ — "uncooked spaghetti" — so that a planar cross-section is the hexagonal disk packing. In general, thick links will include solid cylindrical tubes around proper arcs — "cooked spaghetti" — as well as solid tori around closed curves — "spaghetti-O's." So perhaps this second problem can also be viewed as the natural sequel to the Kepler conjecture (part of Hilbert's 18th problem, now a celebrated theorem of Tom Hales [**6**]) on the density of packing $\mathbb{R}^3$ with equal spheres — "meatballs"!

## 1. Thickness results and conjectures in the round three-sphere

The results in this section were motivated in part by an attempt to better understand the various notions of conformal modulus of a solid torus in $\mathbb{S}^3$, and also the conjecture of Freedman-He-Wang on the Möbius energy minimizer for the Hopf link (see [**1**] and references therein). As it stands, these are preliminary results on thickness, although a goal is to relate thickness and conformal modulus, ultimately getting sharper bounds on the latter.

In $\mathbb{S}^3$ we have this elementary thickness estimate for a single curve:

PROPOSITION 1.1. *The thickness of a tube around a simple closed curve in the round three-sphere cannot exceed that of a maximal tube around a great circle.*

PROOF. Let $\mathcal{K} \subset \mathbb{S}^3$ be the core curve, and observe that the thickness of any osculating circle to $\mathcal{K}$ gives an upper bound for the thickness of $\mathcal{K}$. But a great circle is the thickest round circle in $\mathbb{S}^3$, with thickness $\frac{\pi}{2}$. □

An alternative argument instead uses great circles as secants. The diameter of $\mathcal{K}$ is realized by a secant arc of length at most $\pi$, and so the thickness can be measured from the endpoints of this arc to be (at most) half this length. Using a variation on this argument, and the existence [**7**] of an essential *quadrisecant*, this estimate can be improved by a factor of 2 for a nontrivial thick knot or link in $\mathbb{S}^3$:

THEOREM 1.2. *The thickness of a nontrivial knot or link in the round three-sphere cannot exceed that of the standard geometric Hopf link, whose core consists of two orthogonal great circles.*

PROOF. Let $\mathcal{K}$ denote the nontrivial knot or link. The stereographic projection of $\mathcal{K}$ has an essential quadrisecant in $\mathbb{R}^3$ which pulls back to a round (not necessarily great) circle in $\mathbb{S}^3$ meeting $\mathcal{K}$ in (at least) four points. The shortest of the corresponding four (essential) circular arcs has length at most $\frac{\pi}{2}$. By the Lemma below, such an arc must exit and re-enter the thick tube around $\mathcal{K}$, so half this length is a bound for the thickness of $\mathcal{K}$. ■

LEMMA 1.3. *A secant arc $\mathcal{S}$ which never exits the thick tube around $\mathcal{K}$ is not essential.*

PROOF. Let $\mathcal{B}$ denote the portion of the thick tube around $\mathcal{K}$ which contains $\mathcal{S}$ and lies between normal (great sphere) disks to $\mathcal{K}$ at the endpoints $p, q$ of $\mathcal{S}$; let $\mathcal{K}' = \mathcal{K} \cap \mathcal{B}$. The foliation of $\mathcal{B}$ by normal disks gives a height function on the simple closed curve $\mathcal{K}' \cup \mathcal{S}$ with exactly two critical points ($p$ and $q$), so it bounds an embedded disk whose interior is disjoint from $\mathcal{K}$. By definition, $\mathcal{S}$ is not essential. □

Since Greg Kuperberg's proof for the existence of a quadrisecant [**7**] is fairly involved, it is nice to also have an elementary argument yielding the result of the Theorem in the special case of a two-component link which is not split:

**PROPOSITION 1.4.** *The thickness of a non-split link in the round three-sphere cannot exceed that of the standard geometric Hopf link.*

**PROOF.** This amounts to showing $\mathcal{K} \cup \mathcal{L} \subset \mathbb{S}^3$ is a split link, assuming the simple closed curves $\mathcal{K}$ and $\mathcal{L}$ have distance strictly greater than $\frac{\pi}{2}$ from each other. It suffices to prove that $\mathcal{K}$ is spanned by a singular disk which is disjoint from $\mathcal{L}$. Pick a point $p \in \mathcal{K}$ so that the antipodal point $-p \notin \mathcal{K}$ (since $\mathcal{K}$ can always be approximated by a curve with such a point, there is no loss of generality assuming this for $\mathcal{K}$ itself). Join $p$ to the other points of $\mathcal{K}$ by a one-parameter family of minimizing great circle arcs in $\mathbb{S}^3$. The union of these arcs defines the disk spanning $\mathcal{K}$. Since the length of each arc is at most $\pi$, no point of the disk is more than $\frac{\pi}{2}$ from a point of $\mathcal{K}$, and thus the disk is disjoint from $\mathcal{L}$. □

Note that the thickest geometric Hopf link has the same thickness ($\frac{\pi}{4}$) as the trivial two-component link consisting of parallels on a great $\mathbb{S}^2 \subset \mathbb{S}^3$ at distance $\frac{\pi}{2}$ from each other, and it is not obvious whether this trivial link might be perturbed to something even thicker. It would be interesting to show no other links can realize this bound.

It is clear that the thickest geometric Hopf link in $\mathbb{S}^3$ is the preimage via Hopf projection of a pair of antipodal points on $\mathbb{S}^2$. Similarly, a configuration of $n$ points on $\mathbb{S}^2$ lifts to an $n$-component geometric Hopf link in $\mathbb{S}^3$. Because Hopf projection $\mathbb{S}^3 \to \mathbb{S}^2$ doubles the distance between fibres (Clifford parallel great circles) in $\mathbb{S}^3$ to the distance between the corresponding points in $\mathbb{S}^2$, the thickness of a geometric Hopf link is half the maximal radius for a uniform disk packing around the configuration of points on $\mathbb{S}^2$.

Let $r_n$ denote this maximal radius for a packing of $n$ disks on $\mathbb{S}^2$. (Except for $n \le 6$ and a few other special cases, neither the value of $r_n$, nor the optimal configuration of $n$ points on $\mathbb{S}^2$, is known — for instance, see [**4**].) It is almost surely true that the corresponding geometric Hopf link in $\mathbb{S}^3$ is the thickest in its isotopy class. More daring is the following:

**CONJECTURE 1.5.** *The thickness of any $n$-component link in the round three-sphere cannot exceed $\frac{1}{2}r_n$, the maximal thickness of the geometric Hopf link of $n$ components.*

Similarly for $(m, 2)$-torus knots and links, such as the trefoil ($m = 3$), it is tempting to conjecture the thickest configuration lies on the surface of a flat Clifford torus in $\mathbb{S}^3$ whose aspect ratio is chosen to maximize the strand separation on the torus itself. The reader is encouraged to work this one out for herself!

## 2. On packing density and "The spaghetti-O's conjecture"

Children playing with pennies soon discover the optimal way to pack equal radius disks in the plane is by centering them at the vertices of the hexagonal lattice.

The density of this hexagonal packing — the ratio of area covered to total area — is $\rho_\infty = \frac{\pi}{\sqrt{12}} = .9069\cdots$, and of course this ratio is independent of the radius of the disks packed.

As noted in the previous section, packing $n$ equal radius disks on the unit sphere $\mathbb{S}^2$ is quite a bit trickier. For example, the maximal density $\rho_n$ and maximal radius $r_n$ depend subtly on $n$, and the optimal packings are generally unknown. Clearly

$$\rho_1 = \rho_2 = 1,$$

and spherical trigonometry (with the help of Gauss-Bonnet) lets one compute the local packing density as a decreasing function of ambient curvature (try this as a hard exercise, or see, for example, [**4**] for other approaches) to deduce:

LEMMA 2.1. *The packing density for (at least three) equal-radius disks on the round two-sphere is less than the packing density of the plane:*

$$\rho_n \leq \rho_\infty \quad (n \geq 3).$$

Similarly, one can work out the packing radius $r_n$ for certain small values of $n$; there are also asymptotic formulas [**4**] for both $\rho_n$ and $r_n$ as $n \to \infty$, though this will not be needed in what follows.

How does this apply to compute the packing density for thick tubes around links? Lifting the $n$-disk packing in $\mathbb{S}^2$ via the Hopf map yields a thick Hopf link packing in $\mathbb{S}^3$ with the same density $\rho_n$. For $n = 1$ or 2, this is clearly the densest possible (recall that $\rho_1 = \rho_2 = 1$); but in general, one doesn't expect the packing density for the thick tubes around an $n$-component link to be less than $\rho_n$, since one could make the same link by repeatedly doubling back with a longer, thinner (and thus more dense) tube. Nevertheless, one always expects:

CONJECTURE 2.2. *A thick link in the round three-sphere with more than two components has packing density strictly less than the density $\rho_\infty = \frac{\pi}{\sqrt{12}}$ of hexagonal disk packing. Similarly, $\rho_\infty$ is an upper bound for the asymptotic density for any thick link packing of Euclidean three-space.*

Here *asymptotic density* means the limit supremum of the density within any family of balls whose radius tends to infinity. This might be called "The spaghetti conjecture" — or, as John Sullivan quipped when it was posed to him, "The spaghetti-O's conjecture"! The remainder of this note will be devoted to describing a few examples which illustrate the sharpness of this conjecture, as well as some related conjectures which might lead to a proof.

## 3. Examples of packing thick links, and arguments for the conjecture

Packing $\mathbb{R}^3$ with infinite solid unit-radius cylinders whose cross section is hexagonal disk packing in the plane will clearly realize the density $\rho_\infty = \frac{\pi}{\sqrt{12}} = .9069\cdots$. It is a remarkable theorem of A Bezdek and W Kuperberg [**2**] that this is the densest possible packing of such solid cylinders. Whether or not this is true for finite solid cylinders remains an open question: this may be why it is tricky to get loose, uncooked spaghetti back into the box, or why the author's old station wagon is brim-full of empty soda cans!

At the other extreme from infinitely long spaghetti, consider packings of *bialys* — that is, unit tubes around a unit circle — which can be surprisingly dense. Note that each bialy has volume $2\pi^2 = 19.739\cdots$, by a theorem of Archimedes, and fits inside a solid cylinder of height and radius 2. These cylinders can then be stacked and hexagonally packed to give a bialy packing periodic under the lattice

$$\mathbb{Z}(4, 0, 0) + \mathbb{Z}(2, 2\sqrt{3}, 0) + \mathbb{Z}(0, 0, 2),$$

whose fundamental domain has volume $16\sqrt{3} = 27.712\cdots$, yielding a density of $.7122\cdots$. But the bialys have rounded edges, so shifting adjacent stacks by one unit vertically allows their cores to come a little closer together than 4. Let $c = 2 + \sqrt{3}$ denote this shifted intercore distance. A simple checkerboard pattern leads to a packing of bialys according to the lattice

$$\mathbb{Z}(c, c, 0) + \mathbb{Z}(c, 0, 1) + \mathbb{Z}(0, 0, 2),$$

with slightly greater volume $2c^2 = 27.856\cdots$ and lower density. But this packing has room to let one shear (and rotate) the checkerboard pattern, so that (with $b = \sqrt{c^2 - 4}$) bialys pack in the lattice

$$\mathbb{Z}(4, 0, 0) + \mathbb{Z}(2, b, 1) + \mathbb{Z}(0, 0, 2),$$

with smaller volume $8b = 25.203\cdots$ and higher density $.7830\cdots$. This is denser than any sphere packing, but not very close to the hexagonal packing density. If instead one enlarges the bialy holes — to about the shape of real spaghetti-O's — enough to thread a noodle through each stack, the density can be improved a tiny bit over each of the preceding examples.

Perhaps more remarkable is that one can approach the density of hexagonal disk packing (from below) by using spaghetti-O's of revolution: pack the upper half-plane with disks hexagonally, then revolve around the boundary axis so the disks sweep out solid tori. By the same theorem of Archimedes mentioned above, each solid torus away from the axis has density $\rho_\infty$ within its hexagonal torus cell, but Archimedes' theorem also shows the density for the near-axis bialys is a bit less (it can be computed as $.8950\cdots$). Opening the bialy holes and threading with a noodle does no better.

So why should the spaghetti-O's conjecture be true? There are at least two strategies to verify it. First, one might try to find a clever way to "slice" or "calibrate" a packing of space by thick links, so that the cross sections are (pleated) planes containing lots of disjoint unit disks, and then apply (a pleated version of) the planar disk-packing theorem.

A second approach might be to try proving an even harder result: although the analogue of the Bezdek–Kuperberg packing density theorem is not known for finite height cylinder packings, this has been conjectured (by J B Wilker, for example: *Problem II*, Intuitive Geometry, Colloq. Math. Soc. Janos Bolyai 48 (1987) 700), even for zero-height cylinders. This might lead to a proof of our spaghetti-O's conjecture in an even stronger form, since — unlike the normal disks to the cores of a thick link packing — there would no longer be any coherence required of these disks. Presumably a precise formulation of this stronger conjecture would require the collection of disks to have some kind of transverse measure (so these cylinders really have "infinitesimal" rather than zero height). Unfortunately, the possible presence of curved tubes (such as bialys) in a thick link packing hinders a simple reduction of the spaghetti-O's conjecture to a density result true only for packing cylinders of any positive height.

## Acknowledgements

This material was presented at the April 2001 Las Vegas AMS Meeting, in a special session organized by Jorge Calvo, Ken Millett and Eric Rawdon. Conceived while visiting Jaigyoung Choe at Seoul National University in January 2001, and gestated on several long airline flights between those events, it also has benefited from or stimulated continuing work with my collaborators Jason Cantarella and John Sullivan. The author thanks all these individuals for their interest and hospitality over the past year. Research support from National Science Foundation grant DMS-00-76085 is gratefully acknowledged. Special thanks are also due Christopher Stark, who pointed out the work of W Kuperberg on packing cylinders.

## References

[1] Arnold V I, Khesin B A: *Topological methods in hydrodynamics.* Springer Applied Math. Sci. 125 (1998).

[2] Bezdek A, Kuperberg W: *Maximum density space packing with congruent circular cylinders of infinite length.* Mathematika 37 (1990) 74–80.

[3] Cantarella J, Kusner R B, Sullivan J M: *On the minimum ropelength of knots and links.* ArXiv eprint math.GT/0103224, to appear in Inventiones Math.

[4] Coxeter H S M: *An upper bound for the number of equal nonoverlapping spheres that can touch another of the same size.* Convexity (edited by V Klee), Amer. Math. Soc. Proceedings of Symposia in Pure Math VII (1963) 53–71.

[5] Diao Y: *The lower bounds of the length of thick knots.* Preprint, September 2001.

[6] Hales T C: *Cannonballs and honeycombs.* Notices Amer. Math. Soc. 47 (2000) 440–449.

[7] Kuperberg G: *Quadrisecants of knots and links.* J. Knot Theory Ramifications 3 (1994) 41–50.

CENTER FOR GEOMETRY, ANALYSIS, NUMERICS & GRAPHICS (GANG), DEPARTMENT OF MATHEMATICS, UNIVERSITY OF MASSACHUSETTS, AMHERST, MA 01003
*E-mail address*: kusner@math.umass.edu